\documentclass[PREPRINT,PDF]{tj} 
\usepackage{layout}
\usepackage{amssymb}
\def\hs{\hskip.7pt}
\def\hh{\hskip.4pt}
\def\hn{\hskip-.4pt}
\def\nh{\hskip-.7pt}
\def\nnh{\hskip-1pt}
\def\ns{N}
\def\ip{\eta}
\def\dl{\delta}
\def\ve{\varepsilon}
\def\es{\varSigma}
\def\sg{\varSigma}

\def\bbR{\mathrm{I\!R}}
\def\pc{\mathcal{T}}
\def\tm{{T\hskip-.3ptM}}
\def\tam{{T^*\!M}}
\def\txk{{T\hskip-2.3pt_x^{\phantom i}\hskip-.6ptK}}
\def\txm{{T\hskip-2.3pt_x^{\phantom i}\hskip-.6ptM}}

\def\tym{{T\hskip-2pt_y^{\phantom i}\hskip-.9ptM}}
\def\taxm{{T^{\hh*}_{\hskip-1.7ptx}\hskip-2.4ptM}}
\def\tyn{{T\hskip-2.3pt_y^{\phantom i}\hskip-.9pt\ns}}
\def\tayn{{T^{\hh*}_{\hskip-1.6pty}\hskip-2.4ptN}}
\def\hyp{\hskip.5pt\vbox
{\hbox{\vrule width3ptheight0.5ptdepth0pt}\vskip2.2pt}\hskip.5pt}

\title{Two\hs-\hskip-1.5ptjets of conformal fields along their zero sets}

\shorttitle{Two\hh-\nh jets of conformal fields}

\articletype{Research Article} 

\author{Andrzej Derdzinski\inst{1}\email{andrzej@math.ohio-state.edu}
       }

\institute{\inst{1}
Department of Mathematics, The Ohio State University, Columbus, OH 43210, USA
          }

\abstract{
{\sf 
The connected components of the zero set of any conformal 
vector field $v$, in a pseudo-Riemannian manifold $(M,g)$ of 
arbitrary signature, are of two types, which may be called 
`essential' and `nonessential'. The former consist of points 
at which $v$ is essential, that is, cannot be turned into 
a Killing field by a local conformal change of the metric. 
In a component of the latter type, points at which $v$ is 
nonessential form a relatively-open dense subset that is 
at the same time a totally umbilical submanifold of $(M,g)$. 
An essential component is always a null totally geodesic 
submanifold of $(M,g)$, and so is the set of those points 
in a nonessential component at which $v$ is essential 
(unless this set, consisting precisely of all the singular 
points of the component, is empty). Both kinds of null 
totally geodesic submanifolds arising here carry a 1-form, 
defined up to multiplications by functions without zeros, 
which satisfies a projective version of the Killing equation. 
The conformal-equivalence type of the 2-jet of $v$ is 
locally constant along the nonessential submanifold of a 
nonessential component, and along an essential component 
on which the distinguished 1-form is nonzero. The 
characteristic polynomial of the 1-jet of $v$ is always locally 
constant along the zero set.
}
}


\begin{document}
\maketitle
{\sf 

\section{Introduction}\label{in}
A vector field $\,v\,$ on a pseu\-\hbox{do\hs-}\hskip0ptRiem\-ann\-i\-an 
manifold $\,(M,g)\,$ of dimension $\,n\ge2\,$ is called {\em con\-for\-mal\/} 
if $\,\pounds_vg\,$ equals a function times $\,g$, that is, if for some 
section $\,A\,$ of $\,\mathfrak{so}\hh(\tm)\,$ and some function 
$\,\phi:M\to\bbR$, 
\begin{equation}\label{tnv}
2\hh\nabla\nh v\,\,=\,\,A\,\,+\,\,\phi\hskip2pt\mathrm{Id}\hs,\hskip22pt
\text{or, in coordinates,}\hskip8ptv_{j,\hh k}^{\phantom.}\,
+\,v_{\hh k,\hs j}^{\phantom.}\,=\,\hs\phi\hs g_{jk}^{\phantom.}\hh.
\end{equation}
The covariant 
derivative $\,\nabla\nh v\,$ is treated here as the bundle morphism 
$\,\tm\to\tm\nh\,$ sending any vector field $\,w\,$ to 
$\,\nabla\nh_{\!w}^{\phantom i}\hn v$, and sections of $\,\mathfrak{so}\hh(\tm)\,$ 
are en\-do\-mor\-phisms of $\,\tm\nh$, skew-ad\-joint at every point; 
clearly, $\,\phi=(2/n)\hskip2pt\mathrm{div}\hskip2ptv$.

\ 

If $\,n\ge3$, such $\,v\,$ is known to be uniquely determined by its $\,2$-jet 
at any given point. Determining how the $\,2$-jet of $\,v$ may vary 
along the zero set $\,Z\,$ of $\,v\,$ is thus an obvious initial step 
towards understanding the dynamics of $\,v\,$ near $\,Z$.

\phantom{xx}

\ 

Theorem~\ref{charp} of this paper, which is an easy consequence of 
some facts proved in \cite{derdzinski}, deals with the $\,1$-jet of $\,v$, 
establishing a restriction on its variability: the characteristic polynomial of 
$\,\nabla\nh v\,$ must be locally constant on $\,Z$.

\ 

A point $\,x\in Z\,$ is called {\em nonessential\/} if some local con\-for\-mal 
change of the metric at $\,x\,$ turns $\,v\,$ into a Kil\-ling field, and {\em 
essential\/} otherwise. A connected components of $\,Z\,$ is either {\em 
essential\/} (if it consists of essential points only) or {\em nonessential\/} 
(when it contains some nonessential points, possibly along with essential ones).

\ 

The next main result, Theorem~\ref{essen}, explores structural properties of 
components of $\,Z$. Every essential component turns out to be a null totally 
geodesic sub\-man\-i\-fold, and so is, when nonempty, the 
pos\-si\-bly-dis\-con\-nect\-ed set $\,\sg\,$ of 
essential points in any given nonessential component $\,\ns$. At the same time, 
$\,\sg\,$ coincides with the set of singular points of $\,\ns$, while 
$\,\ns\smallsetminus\sg\,$ is a totally umbilical sub\-man\-i\-fold. The 
tangent spaces of these sub\-man\-i\-folds at all points $\,x\,$ are explicitly 
described in terms of $\,\nabla\nh v_x^{\phantom i}$ and $\,d\phi_x^{\phantom i}$.

\ 

For $\,\ns\,$ and $\,\sg\,$ as above, let the same symbol $\,\es\,$ also stand 
for an essential component of $\,Z$. Section~\ref{is} discusses geometric 
structures on $\,\ns\smallsetminus\sg\,$ and both types of $\,\es$, naturally 
induced by the underlying con\-for\-mal structure of $\,(M,g)$. They consist 
of a con\-stant-rank, pos\-si\-bly-de\-gen\-er\-ate con\-for\-mal structure on 
$\,\ns\smallsetminus\sg\,$ along with its null\-space distribution, a 
projective structure on $\,\es$, and a $\,1$-form $\,\xi\,$ on $\,\es\,$ 
defined only up to multiplications by functions without zeros. Their basic 
properties are listed in Proposition~\ref{pties}.

\ 

Finally, Section~\ref{tj} addresses the question, mentioned above, of 
variability of the $\,2$-jet of $\,v\,$ along $\,Z$. The 
con\-for\-mal-e\-quiv\-a\-lence type of the $\,2$-jet is proved to 
be locally constant in $\,\ns\smallsetminus\sg\,$ and, generically, in 
$\,\es$. The word `generically' means here {\em in any component of\/ 
$\,\es\,$ on which\/ $\,\xi\,$ is not identically zero}. Examples show that, 
in the case of $\,\es$, some form of the `generic' assumption is necessary. 
On the other hand, if $\,\sg\subset\ns\,$ is nonempty, the equivalence types 
at points of $\,\sg\,$ are always different from those realized in 
$\,\ns\smallsetminus\sg$.

\section{Preliminaries}\label{pr}
Manifolds {\em need not be connected}. However, their connected components 
must all have the same dimension. Sub\-man\-i\-fold are always endowed with 
the subset topology. All manifolds, mappings, bundles and their sections, 
including tensor fields and functions, are of class $\,C^\infty\nnh$. The 
symbol $\,\nabla\,$ denotes both the Le\-vi-Ci\-vi\-ta connection of a given 
pseu\-\hbox{do\hs-}\hskip0ptRiem\-ann\-i\-an metric $\,g\,$ on a manifold 
$\,M$, and the $\,g$-grad\-i\-ent. Thus, for a vector field $\,u\,$ and a 
function $\,\tau\,$ on $\,M$, we have $\,d_u\tau=g(u,\nabla\tau)$.

\ 

Given a sub\-man\-i\-fold $\,K\,$ of a manifold $\,M$, we denote by 
$\,T\hskip-2.7pt_K^{\phantom i}\nnh M\,$ the restriction of $\,\tm\,$ to 
$\,K$. The {\em normal bundle\/} of $\,K\,$ is defined, as usual, to be the 
quotient vector bundle $\,T\hskip-2.7pt_K^{\phantom i}\nnh M\nnh/T\nh K$. Any 
fixed tor\-sion-free connection $\,\nabla\,$ on $\,M\,$ gives rise to the 
{\em second fundamental form\/} of $\,K$, which is a section $\,b\,$ of 
$\,[\tam]^{\odot2}\nnh\otimes T\hskip-2.7pt_K^{\phantom i}\nnh M\nnh/T\nh K\,$ 
(so that, at every $\,x\in K$, the mapping $\,b_x:\txk\times\txk\to\txm\nnh/\txk\,$ 
is bi\-lin\-e\-ar and symmetric). We have 
\begin{equation}\label{bxw}
b(\dot x,w)\,\,=\,\,\pi\nabla_{\!\dot x}^{\phantom i}w
\end{equation}
whenever $\,t\mapsto w(t)$ 
is a vector field tangent to $\,K\,$ along a curve $\,t\mapsto x(t)\,$ in 
$\,K$, with $\,\pi:\tm\to T\hskip-2.7pt_K^{\phantom i}\nnh M\nnh/T\nh K\,$ 
denoting the quotient projection. When $\,b=0\,$ identically, $\,K\,$ is 
said to be {\em totally geodesic\/} relative to $\,\nabla\nnh$. If $\,\nabla\,$ is 
the Le\-vi-Ci\-vi\-ta connection of a 
pseu\-\hbox{do\hs-}\hskip0ptRiem\-ann\-i\-an metric $\,g\,$ and 
$\,b=g_{\nh K}^{\phantom i}\nh\otimes u\,$ for some section $\,u\,$ of 
$\,T\hskip-2.7pt_K^{\phantom i}\nnh M\nnh/T\nh K$, where 
$\,g_{\nh K}^{\phantom i}$ is the restriction of $\,g$ to $\,K$, one calls 
$\,K\,$ {\em totally umbilical\/} in $\,(M,g)$. This last property is 
con\-for\-mal\-ly invariant: changing $\,g\,$ to $\,e^\tau\hskip-2ptg\,$ 
causes $\,b\,$ to be replaced by 
$\,b\hskip1pt-g_{\nh K}^{\phantom i}\nh\otimes\pi\nabla\nh\tau\nh/2$. 

\ 

As shown by Weyl \cite[p.~100]{weyl}, two 
tor\-sion\-free connections on a manifold $\,M\,$ are {\em projectively 
equivalent}, in the sense of having the same re-pa\-ram\-e\-trized 
geodesics, if and only if their difference $\,E\,$ can be written as 
$\,E\hn=\hs\theta\hn\odot\nh\text{\rm Id}\,$ for some $\,1$-form $\,\theta\,$ 
on $\,M\,$ (in coordinates: $\,2\hh E_{jk}^l
=\theta_j^{\phantom i}\delta_k^l\nh+\theta_k^{\phantom i}\delta_j^l$). On the 
other hand, given a pseu\-\hbox{do\hs-}\hskip0ptRiem\-ann\-i\-an metric 
$\,g$ on $\,M$, with the Le\-vi-Ci\-vi\-ta connection $\,\nabla\nnh$, 
and a function $\,\tau:M\to\bbR$, the con\-for\-mal\-ly related metric 
$\,e^\tau\hskip-2ptg\,$ has the Le\-vi-Ci\-vi\-ta connection 
$\,\nabla\nh+E$, where $\,E=d\tau\odot\text{\rm Id}-g\otimes\nabla\tau\nh/2$. 
Thus, for any null totally geodesic sub\-man\-i\-fold $\,\es\,$ of $\,M$, 
the connections on $\,\es\,$ induced by the Le\-vi-Ci\-vi\-ta connections 
of $\,g\,$ and $\,e^\tau\hskip-2ptg\,$ are projectively equivalent.

\ 

For every con\-for\-mal vector field $\,v\,$ on a 
pseu\-\hbox{do\hs-}\hskip0ptRiem\-ann\-i\-an manifold $\,(M,g)\,$ of dimension 
$\,n\ge3\,$ and any vector fields $\,u,w\,$ on $\,M\,$ one has the well-known 
equalities of bundle mor\-phisms $\,\tm\to\tm\,$ and functions $\,M\to\bbR\hs$:
\begin{equation}\label{nwn}
\begin{array}{l}
2\hh\nabla_{\!u}^{\phantom i}\nabla\nh v
=2\hs R(v\wedge u)+d\phi\nh\otimes\nh u
-g(u,\,\cdot\,)\nnh\otimes\nnh\nabla\nh\phi
+g(u,\nabla\nh\phi)\hs\mathrm{Id}\hs,\\
(1-n/2)[\nabla d\phi](u,w)
=S(u,\nabla\nh_{\!w}^{\phantom i}\hn v)+S(w,\nabla_{\!u}v)
+[\hh\nabla_{\!v}S\hh](u,w)\hs,
\end{array}
\end{equation}
cf.\ \cite[formula (22))]{derdzinski}, where $\,R\,$ and $\,S\,$ are the 
curvature and Schouten tensors. In coordinates, (\ref{nwn}) reads 
$\,2\hh v^{\hs l}{}_{,\hh kj}\nh=2\hh R_{\hh pjk}{}^{\hh l}v^{\hh p}\nh
+\phi_{\nh,\hh k}^{\phantom i}\delta_{\nh j}^{\hs l}\nh
-\phi^{\hs,\hh l}g_{jk}^{\phantom.}\nh
+\phi_{\nh,\hs j}^{\phantom i}\delta_k^{\hs l}\,$ and 
$\,\hs(1-n/2)\phi_{\nh,\hs jk}\nh=S_{jp}v^{\hs p}{}_{,\hh k}\nh
+S_{kp}v^{\hs p}{}_{,\hs j}\nh+S_{jk,\hh p}v^{\hh p}\nnh$.
\begin{remark}\label{phcnf}For a con\-for\-mal vector field $\,v\,$ on a 
pseu\-\hbox{do\hs-}\hskip0ptRiem\-ann\-i\-an manifold $\,(M,g)\,$ and any 
function $\,\tau:M\to\bbR$, the con\-for\-mal\-ly equivalent metric 
$\,e^\tau\hskip-2ptg\,$ satisfies, along with $\,v$, the analog of (\ref{tnv}) 
in which the role of $\,\phi\,$ is played by $\,\phi+d_v\tau$. In fact, 
(\ref{tnv}) is equivalent to $\,\pounds_vg=\phi\hs g$, while 
$\,\pounds_v(e^\tau\hskip-2ptg)=e^\tau\nnh\pounds_vg
+(d_v\tau)\hs e^\tau\hskip-2ptg$. At a point $\,x\,$ such that $\,v_x\nh=0$, 
switching from $\,g\,$ to $\,e^\tau\hskip-2ptg\,$ thus results in replacing 
$\,d\phi_x$ by $\,d\phi_x\nh+(d\hh\tau_x)\nabla\nh v_x$. 
\end{remark}
\begin{remark}\label{kilng}A Kil\-ling field $\,v\,$ and any vector field 
$\,u\,$ on a pseu\-\hbox{do\hs-}\hskip0ptRiem\-ann\-i\-an manifold satisfy 
(\ref{nwn}) with $\,\phi=0$. Thus, $\,\nabla\nh v\,$ is parallel along any 
curve to which $\,v\,$ is tangent, such as an integral curve of $\,v\,$ or a 
curve of zeros of $\,v$.
\end{remark}

\section{The zero set $\,Z\,$ of a conformal field $\,v$}\label{zs}
In addition to the function 
$\,\phi=(2/n)\hskip2pt\mathrm{div}\hskip2ptv:M\to\bbR\,$ appearing in 
(\ref{tnv}), let us also consider
\begin{equation}\label{ziz}
{}\mathsf{the\ zero\ set\ }\,Z\,\mathsf{\ of\ a\ con\-for\-mal\ vector\ field\ }\,v\,
\mathsf{\ on\ a\ pseu\-\hbox{do\hs-}\hskip0ptRiem\-ann\-i\-an\ manifold\ }\,(M,g),
\hskip6pt\dim M=n\ge3\hh.
\end{equation}
If $\,x\in Z$, the {\em simultaneous kernel at\/} $\,x\,$ of the differential $\,d\phi\,$ 
and the bundle morphism $\,\nabla\nh v:\tm\to\tm\,$ is the space
\begin{equation}\label{hxe}
\mathcal{H}_x^{\phantom i}\,\,
=\,\,\hs\mathrm{Ker}\hskip1.7pt\nabla\nh v_x^{\phantom i}
\cap\hs\mathrm{Ker}\hskip2.2ptd\phi_x^{\phantom i}\hh.
\end{equation}
When $\,x\,$ is fixed, the symbol $\,H\,$ may be used instead of 
$\,\mathcal{H}_x^{\phantom i}$.

\ 

As in 
\cite{belgun-moroianu-ornea}, we call $\,x\in Z\,$ a {\em nonessential\/} zero 
of $\,v\,$ if $\,v\,$ restricted to a suitable neighborhood of $\,x\,$ is a 
Kil\-ling field for some metric con\-for\-mal to $\,g$. When no such 
neighborhood and metric exist, the zero of $\,v\,$ at $\,x\,$ is said to be 
{\em essential}.

\ 

By a {\em nonsingular point\/} of $\,Z\,$ we mean any $\,x\in Z\,$ such 
that, for some neighborhood $\,\,U\,$ of $\,x\,$ in $\,M$, the intersection 
$\,Z\cap U\,$ is a sub\-man\-i\-fold of $\,M$. Points of $\,Z\,$ not having a 
neighborhood with this property will be called {\em singular}. 

\ 

For $\,(M,g),\hs v,\hh Z\,$ as above, a point $\,x\in Z$, and the 
exponential mapping $\,\exp\hskip-2pt_x^{\phantom i}$ of $\,g\,$ 
at $\,x$, we will repeadtedly consider
\begin{equation}\label{suf}
\begin{array}{l}
\mathsf{any\ \hskip1ptsufficiently\ \hskip1ptsmall\ \hskip1ptneighborhoods\ \hskip1pt}\,\,U\,\mathsf{\ \hskip1ptof\ \hskip1pt}\,0\,
\mathsf{\ \hskip1ptin\ \hskip1pt}\,\txm\,\mathsf{\ \hskip1ptand\ \hskip1pt}\,\,U'\mathsf{\ \hskip1ptof\ \hskip1pt}\,x\,\mathsf{\ \hskip1pt
in\ \hskip1pt}\,M\,\mathsf{\ \hskip1ptsuch\ \hskip1ptthat}\\
U\,\mathsf{\ is\ a\ union\ of\ line\ segments\ emanating\ from\ }\,0\,\mathsf{\ 
and\ }\,\exp\hskip-2pt_x^{\phantom i}\mathsf{\ is\ a\ dif\-feo\-mor\-phism\ }
\,\,U\to U'.
\end{array}
\end{equation}
\begin{theorem}[Kobayashi \cite{kobayashi}]\label{zrklf}For\/ 
$\,(M,g),\hs v,\hh Z,\,U,U'$ and\/ $\,H=\mathcal{H}_x^{\phantom i}$ as in\/ 
{\rm(\ref{ziz})} -- {\rm(\ref{suf})}, let\/ $\,v\,$ also be a Kil\-ling field. 
Then
\[ 
Z\cap U'\,\,=\,\,\hs\exp\hskip-2pt_x^{\phantom i}[H\cap U\hs]\hs,\hskip24pt
\text{with}\hskip12ptH=\hs\mathrm{Ker}\hskip1.7pt\nabla\nh v_x^{\phantom i}
\hskip8pt\text{since\ }\,\,\,\phi=0\,\,\,\text{in\/ {\rm(\ref{tnv})}}.
\]
Thus, the connected components of\/ $\,Z\,$ are totally geodesic 
sub\-man\-i\-folds of even co\-di\-men\-sions in $\,(M,g)$.
\end{theorem}
\begin{theorem}[Beig \cite{beig,capocci}]\label{esszr}Let\/ $\,Z\,$ be the 
zero set of a con\-for\-mal vector field\/ $\,v\,$ on a 
pseu\-\hbox{do\hs-}\hskip0ptRiem\-ann\-i\-an manifold\/ $\,(M,g)\,$ of 
dimension\/ $\,n\ge3$. A point\/ $\,x\in Z\,$ is nonessential if and only if
\begin{equation}\label{nes}
\phi(x)=0\,\,\,\text{\ and\ }\,\,\,
\nabla\nh\phi_x^{\phantom i}\in\nabla\nh v_x^{\phantom i}(\txm)\hh.
\end{equation}
for the function\/ $\,\phi=(2/n)\hskip2pt\mathrm{div}\hskip2ptv:M\to\bbR\,$ 
appearing in\/ {\rm(\ref{tnv})}. In other words, $\,x\in Z\,$ is essential if 
and only if
\begin{equation}\label{ess}
{}\text{either}\hskip6pt\phi(x)\hs\ne\hs0\hs,\hskip9pt\text{or}\hskip8pt
\phi(x)\hs=\hs0\hskip6pt\text{and}\hskip5pt
\nabla\nh\phi_x^{\phantom i}\hs\notin\hs\nabla\nh v_x^{\phantom i}(\txm)\hs.
\end{equation}
\end{theorem}
\begin{theorem}[Derdzinski \cite{derdzinski}]\label{zrset}
Suppose that\/ $\,v\,$ is a con\-for\-mal vector field on a 
pseu\-\hbox{do\hs-}\hskip0ptRiem\-ann\-i\-an manifold\/ $\,(M,g)\,$ of 
dimension\/ $\,n\ge3$. If\/ $\,Z$ is the zero set of\/ $\,v$, while\/ 
$\,x\in Z\,$ satisfies\/ {\rm(\ref{ess})}, and\/ 
$\,C=\{w\in\txm:g_x^{\phantom i}(w,w)=0\hs\}\,$ stands for the null cone, 
then, with\/ $\,\,U,U'$ as in\/ {\rm(\ref{suf})} and\/ 
$\,H=\hs\mathrm{Ker}\hskip1.7pt\nabla\nh v_x^{\phantom i}
\cap\hs\mathrm{Ker}\hskip2.2ptd\phi_x^{\phantom i}$,
\begin{equation}\label{zcu}
Z\cap U'\,\,=\,\,\hs\exp\hskip-2pt_x^{\phantom i}[\hs C\cap H\cap U\hs]\hs.
\end{equation}
In addition, $\,\phi=(2/n)\hskip2pt\mathrm{div}\hskip2ptv\,$ is constant along 
each connected component of\/ $\,Z$.
\end{theorem}
Given a con\-for\-mal vector field $\,v\,$ on a 
pseu\-\hbox{do\hs-}\hskip0ptRiem\-ann\-i\-an manifold $\,(M,g)\,$ of 
dimension $\,n\ge3$, and a parallel vector field 
$\,t\mapsto w(t)\in T\hskip-2.3pt_{y(t)}\hskip-.9ptM\,$ along a geodesic 
$\,t\mapsto y(t)\,$ contained in the zero set $\,Z\,$ of $\,v$, such that 
$\,g(\dot y,w)=0$, replacing $\,u\,$ in (\ref{nwn}) with $\,\dot y\,$ we 
obtain
\begin{equation}\label{nyw}
\mathsf{(a)}\hskip9pt2\hh\nabla_{\!\dot y}\nabla\nh_{\!w}^{\phantom i}\hn v
=g(w,\nabla\nh\phi)\hs\dot y\hs,\hskip18pt
\mathsf{(b)}\hskip9pt(1-n/2)[g(w,\nabla\nh\phi)]\hs\dot{\,}
=S(\dot y,\nabla\nh_{\!w}^{\phantom i}\hn v)\hs.
\end{equation}
(The other terms vanish since $\,v=\nabla_{\!\dot y}v=0\,$ at $\,y(t)$, while  
$\,g(\dot y,\nabla\nh\phi)=0\,$ due to the final clause of 
Theorem~\ref{zrset}.)
\begin{remark}\label{lpcon}
In view of Theorems~\ref{zrklf} --~\ref{zrset}, $\,Z\,$ in (\ref{ziz}) is 
always locally pathwise connected. Thus, the connected components of $\,Z$ 
are pathwise connected, closed subsets of $\,M$.
\end{remark}
\begin{remark}\label{tlumb}
Away from singularities, the connected components of $\,Z\,$ are totally 
umbilical sub\-man\-i\-folds of $\,(M,g)$, and their co\-di\-men\-sions are 
even unless the component is a null totally geodesic sub\-man\-i\-fold.

\ 

In fact, for the connected components of the set of nonsingular zeros of 
$\,v$, this easily follows from Theorems~\ref{zrklf} --~\ref{zrset}; see also 
\cite[Theorem 1.1]{derdzinski}. That the connected components of $\,Z$, with 
the singularities removed, are sub\-man\-i\-folds as well (in other words, 
their own connected components all have the same dimension) is also immediate 
from Theorems~\ref{zrklf} --~\ref{zrset}: by (\ref{zcu}), the set of singular 
points in $\,Z\cap U'$ coincides with 
$\,\exp\hskip-2pt_x^{\phantom i}[\hs H\cap H^\perp\hskip-1.8pt\cap U\hs]\,$ 
when (\ref{ess}) holds and the metric $\,g_x^{\phantom i}$ restricted to 
$\,H\,$ is not sem\-i\-def\-i\-nite, and is empty otherwise, while, in the 
former case, all components of $\,(C\smallsetminus H^\perp)\cap H\,$ are 
clearly of dimension $\,\dim H-1$. 
\end{remark}

\section{The characteristic polynomial of $\,\nabla\nh v$}\label{cp}
Given a tor\-sion-free connection $\,\nabla\,$ on an 
$\,n\hh$-\hskip.7ptdi\-men\-sion\-al manifold $\,M$, and a vector field 
$\,v\,$ on $\,M$, we denote by $\,\mathcal{P}\hskip-3pt_n^{\phantom j}$ the 
space of real all polynomials in one variable of degrees not exceeding $\,n$, 
and by $\,\chi(\nabla\nh v)\,$ the function 
$\,M\to\mathcal{P}\hskip-3pt_n^{\phantom j}$ 
assigning to each $\,x\in M\,$ the characteristic polynomial of the 
en\-do\-mor\-phism $\,\nabla\nh v_x^{\phantom i}:\txm\to\txm$.
\begin{lemma}[Derdzinski \cite{derdzinski}, Lemma 12.2(b)--(iii)]\label{ctchp}
If a con\-for\-mal vector field\/ $\,v\hs$ on a 
pseu\-\hbox{do\hs-}\hskip0ptRiem\-ann\-i\-an manifold\/ $\,(M,g)\,$ is tangent 
to a null geodesic segment\/ $\,\varGamma\nnh$, and\/ $\,\phi\,$ appearing 
in\/ {\rm(\ref{tnv})} is constant along\/ $\,\varGamma\nnh$, then\/ 
$\,\chi(\nabla\nh v)\,$ is constant along\/ $\,\varGamma$ as well.
\end{lemma}
\begin{theorem}\label{charp}Let\/ $\,Z\,$ be the zero set of a con\-for\-mal 
vector field\/ $\,v\,$ on a pseu\-\hbox{do\hs-}\hskip0ptRiem\-ann\-i\-an 
manifold\/ $\,(M,g)\,$ of dimension\/ $\,n\ge3$. Then\/ 
$\,\chi(\nabla\nh v):M\to\mathcal{P}\hskip-3pt_n^{\phantom j}$ is constant on 
every connected component of\/ $\,Z\,$ and, consequently, so is\/ 
$\,\phi=(2/n)\hskip2pt\mathrm{div}\hskip2ptv:M\to\bbR$.
\end{theorem}
\begin{proof}
{\sf 
We fix $\,x\in Z\,$ and show that $\,\chi(\nabla\nh v)$, at zeros 
of $\,v\,$ near $\,x$, is the same as at $\,x$, cf.\ Remark~\ref{lpcon}.

\ 

First, if $\,x\,$ is a nonessential zero of $\,v$, changing the metric 
con\-for\-mal\-ly near $\,x$, we may assume that $\,v\,$ is a Kil\-ling field. 
By Theorem~\ref{zrklf}, the nearby zeros of $\,v\,$ then form a 
sub\-man\-i\-fold $\,K\,$ of $\,M$, while, according to Remark~\ref{kilng}, 
$\,\nabla\nh v$ is parallel along $\,K$. This proves our assertion for 
nonessential zeros $\,x$.

\ 

Now let the zero of $\,v\,$ at $\,x\,$ be essential. Theorem~\ref{esszr} then 
gives (\ref{ess}). In view of Theorem~\ref{zrset}, every nearby point of 
$\,Z\,$ is joined to $\,x\,$ by a null geodesic segment $\,\varGamma\hs$ 
contained in $\,Z$. Our claim about $\,\phi\,$ follows in turn from the final 
clause of Theorem~\ref{zrset}. Constancy of $\,\chi(\nabla\nh v)\,$ along 
$\,\varGamma\hs$ is therefore immediate from Lemma~\ref{ctchp}.
}
\end{proof}

\section{Essential and nonessential components of $\,Z$}\label{en}
By the {\em components\/} of the set $\,Z\,$ appearing in (\ref{ziz}) we mean 
its (pathwise) connected components, cf.\ Remark~\ref{lpcon}.

\ 

A component of $\,Z\,$ will be called {\em essential\/} if all of its 
points are essential zeros of $\,v$, as defined in Section~\ref{zs}. 
Otherwise, the component is said to be {\em nonessential}.

\ 

This definition allows a nonessential component $\,\ns\hs$ to contain 
some essential zeros of $\,v$. However, as shown in Theorem~\ref{essen}(v) 
below, the set of nonessential zeros in $\,\ns\hs$ is relatively open and 
dense.

\ 

Unlike the components, sub\-man\-i\-folds---such as $\,\es\,$ and 
$\,\ns\smallsetminus\sg\,$ in case (b) of Theorem~\ref{essen}---need not be 
connected, although all of their components are required to have the same 
dimension.
\begin{theorem}\label{essen}Given a con\-for\-mal vector field\/ $\,v\,$ on a 
pseu\-\hbox{do\hs-}\hskip0ptRiem\-ann\-i\-an manifold\/ $\,(M,g)\,$ of 
dimension\/ $\,n\ge3$, suppose that
\begin{enumerate}
  \def\theenumi{{\rm\alph{enumi}}}
\item[{\rm(a)}] $\es\,$ is an essential component of\/ $\,Z$, or
\item[{\rm(b)}] $\es\,$ is the set of essential points in a nonessential 
component $\,\ns\hs$ of\/ $\,Z$,
\end{enumerate}
where\/ $\,Z\,$ is the zero set of\/ $\,v$. Then, with\/ 
$\,\mathcal{H}_x^{\phantom i}\,
=\,\hs\mathrm{Ker}\hskip1.7pt\nabla\nh v_x^{\phantom i}
\cap\mathrm{Ker}\hskip2.2ptd\phi_x^{\phantom i}$,
\begin{enumerate}
  \def\theenumi{{\rm\roman{enumi}}}
\item[{\rm(i)}] $\sg$, if nonempty, is a null totally geodesic 
sub\-man\-i\-fold of $\,(M,g)$, closed as a subset of $\,M\nh$,
\item[{\rm(ii)}] $T\hskip-2.3pt_x^{\phantom i}\es\,
=\,\hs\mathcal{H}_x^{\phantom i}\nh\cap\mathcal{H}_x^\perp\,$ at every point\/ 
$\,x\,$ of\/ $\,\es\nh$,
\item[{\rm(iii)}] for any\/ $\,x\in\es\,$ the metric\/ $\,g_x^{\phantom i}$ 
restricted to\/ $\,\mathcal{H}_x^{\phantom i}$ is sem\-i\-def\-i\-nite in 
case\/ {\rm(a)}, non-sem\-i\-def\-i\-nite in case\/ {\rm(b)}.
\end{enumerate}
Finally, in case\/ {\rm(b)}, with singular points and\/ $\,C,H\,$ defined as 
in Section\/~{\rm\ref{zs}} and Theorem\/~{\rm\ref{zrset}},
\begin{enumerate}
  \def\theenumi{{\rm\roman{enumi}}}
\item[{\rm(iv)}] $\sg\,$ consists of singular, $\ns\smallsetminus\sg\,$ of 
nonsingular zeros of\/ $\,v\,$ in\/ $\,\ns\nnh$,
\item[{\rm(v)}] $\ns\smallsetminus\sg\,$ is a totally umbilical 
sub\-man\-i\-fold of $\,M$, while the sign pattern of\/ $\,g\,$ restricted 
to\/ $\,\ns\smallsetminus\sg$, including its rank\/ $\,r$, is the same at all 
points, and\/ $\,\dim\hs(\ns\smallsetminus\sg)\hs-\hs\dim\hs\sg\hs=\hs r+1$,
\item[{\rm(vi)}] whenever\/ $\,y\in\ns\smallsetminus\sg\,$ one has\/ 
$\,T\hskip-2.3pt_y^{\phantom i}(\ns\smallsetminus\sg)
=\hs\mathrm{Ker}\hskip1.7pt\nabla\nh v_y^{\phantom i}$, and\/ 
$\,\mathrm{rank}\hskip1.7pt\nabla\nh v_y^{\phantom i}\hs
=\hs2\hs+\hs\mathrm{rank}\hskip1.7pt\nabla\nh v_x^{\phantom i}\,$ if\/ 
$\,x\in\sg$,
\item[{\rm(vii)}] for\/ $\,x\in\sg\,$ and sufficiently small $\,\,U,U'$ 
in\/ {\rm(\ref{suf})}, $\,\sg\cap U'\nh
=\exp\hskip-2pt_x^{\phantom i}[\hs H\cap H^\perp\hskip-1.8pt\cap U\hs]\,$ 
and\/ $\,\ns\cap U'\nh=\exp\hskip-2pt_x^{\phantom i}[\hs C\cap H\cap U\hs]$.
\end{enumerate}
\end{theorem}
\begin{proof}
{\sf 
As a consequence of Theorem~\ref{esszr}, for $\,x\in Z\,$ and 
$\,H=\mathcal{H}_x^{\phantom i}$  
there are three possibilities:
\begin{enumerate}
  \def\theenumi{{\rm\Alph{enumi}}}
\item[{\rm($\alpha$)}] $x\,$ is a nonessential zero of $\,v$, that is, 
(\ref{nes}) holds,
\item[{\rm($\beta$)}] $x\,$ is essential and the metric $\,g_x^{\phantom i}$ 
is sem\-i\-def\-i\-nite on $\,H$,
\item[{\rm($\gamma$)}] $x\,$ is essential, $\,g_x^{\phantom i}$ restricted to 
$\,H$ is not 
sem\-i\-def\-i\-nite, $\,\phi(x)=0\,$ and 
$\,\nabla\nh\phi_x^{\phantom i}\hs\notin\hs\nabla\nh v_x^{\phantom i}(\txm)$.
\end{enumerate}
The assertions about $\,\phi=(2/n)\hskip2pt\mathrm{div}\hskip2ptv\,$ in 
($\gamma$) follow from (\ref{ess}); note that, if $\,\phi(x)\,$ were nonzero, 
$\,\mathrm{Ker}\hskip1.7pt\nabla\nh v_x^{\phantom i}$ would be a null 
sub\-space of $\,\txm\,$ (as an eigen\-space of the skew-ad\-joint 
en\-do\-mor\-phism $\,A_x^{\phantom i}$ in (\ref{tnv}) for a nonzero 
eigen\-value), and so $\,g_x^{\phantom i}$ would be sem\-i\-def\-i\-nite on 
$\,H\subset\mathrm{Ker}\hskip1.7pt\nabla\nh v_x^{\phantom i}$. In view 
of Theorem~\ref{zrklf} and \cite[second paragraph on p.\ 22]{derdzinski}, 
\begin{equation}\label{sng}
x\,\,\mathsf{\ is\ nonsingular\ in\ cases\ }(\alpha)\mathsf{\ and\ 
}(\beta),\mathsf{\ but\ singular\ in\ case\ }(\gamma).
\end{equation}
If $\,Z\cap U'$ satisfy ($\gamma$) then, for $\,\,U,U'$ as in (\ref{suf}), 
$\,\sg\hh'\nh
=\exp\hskip-2pt_x^{\phantom i}[\hs H\cap H^\perp\hskip-1.8pt\cap U\hs]\,$ and 
$\,\ns'\nh
=\exp\hskip-2pt_x^{\phantom i}[\hs C\cap H\cap U\hs]\smallsetminus\sg\hh'$ are 
sub\-man\-i\-folds of $\,M\,$ such that
\begin{equation}\label{tyn}
\tyn'\,=\,\,\mathrm{Ker}\hskip1.7pt\nabla\nh v_y^{\phantom i}\hskip7pt
\mathsf{and}\hskip9pt\mathrm{rank}\hskip1.7pt\nabla\nh v_y^{\phantom i}\,
=\,\hs2\hs\,+\,\mathrm{rank}\hskip1.7pt\nabla\nh v_x^{\phantom i}\hskip7pt\mathsf{for\ 
every}\hskip6pty\in\ns'.
\end{equation}
In fact, $\,H\cap H^\perp$ clearly is the set of singular points in 
$\,C\cap H$. (For more details, see \cite[Remark 6.2(a)]{derdzinski}.) To 
verify (\ref{tyn}) for sufficiently small $\,\,U,U'\nnh$, note that 
$\,\tyn'\nh\subset\hs\mathrm{Ker}\hskip1.7pt\nabla\nh v_y^{\phantom i}$ as 
$\,\ns'\nh\subset Z$, while $\,\dim\ns'\nh=\dim H-1
=\dim\hs\mathrm{Ker}\hskip1.7pt\nabla\nh v_x^{\phantom i}-2\,$ due to the 
definiton of $\,H=\mathcal{H}_x^{\phantom i}$ and the last relation in 
($\gamma$). This yields 
$\,\dim\hs\mathrm{Ker}\hskip1.7pt\nabla\nh v_x^{\phantom i}-2
\le\dim\hs\mathrm{Ker}\hskip1.7pt\nabla\nh v_y^{\phantom i}$ or, 
equivalently, 
$\,\mathrm{rank}\hskip1.7pt\nabla\nh v_x^{\phantom i}
\le\mathrm{rank}\hskip1.7pt\nabla\nh v_y^{\phantom i}
\le\,2\hs+\mathrm{rank}\hskip1.7pt\nabla\nh v_x^{\phantom i}$ (where, 
for $\,y\,$ near $\,x$, we have also used sem\-i\-con\-ti\-nu\-i\-ty 
of the rank). The two inequalities cannot be both strict, as both ranks 
are even: $\,\nabla\nh v_x^{\phantom i}$ and $\,\nabla\nh v_y^{\phantom i}$ 
are skew-adjoint by (\ref{tnv}), with
\begin{equation}\label{fez}
\phi\,=\,0\hskip8pt\mathsf{\ on}\hskip6pt\ns'\nh\cup\sg\hh'
\end{equation}
in view the final clause of Theorem~\ref{zrset} and ($\gamma$). If we now 
did not have $\,\mathrm{rank}\hskip1.7pt\nabla\nh v_y^{\phantom i}
=\hs2\hs+\mathrm{rank}\hskip1.7pt\nabla\nh v_x^{\phantom i}$ for all 
$\,y\in\ns'$ close to $\,x$, there would be a sequence 
of points $\,y\in\ns'$ such that 
$\,\mathrm{rank}\hskip1.7pt\nabla\nh v_y^{\phantom i}
=\mathrm{rank}\hskip1.7pt\nabla\nh v_x^{\phantom i}$, converging to $\,x\,$ 
and, by continuity, ($\gamma$) with $\,x\,$ replaced by $\,y\,$ would hold 
for all but finitely many of its terms $\,y$. (They would be essential in 
view of (\ref{ess}), with $\,\phi(y)=0\,$ by (\ref{fez}).) The result would be 
a contradiction, as $\,y\,$ would then be singular by (\ref{sng}), yet at the 
same time nonsingular since, in view of Theorem~\ref{zrset} the 
sub\-man\-i\-fold $\,\ns'$ of $\,M$, containing $\,y$, is a relatively open 
subset of $\,Z$.

\ 

Furthermore, for $\,\sg\hh'\nnh,\ns'$ chosen as above, with 
sufficiently small $\,\,U,U'\nnh$,
\begin{equation}\label{alg}
\mathsf{points\ of\ }\,\sg\hh'\mathsf{\ have\ property\ }(\gamma),
\mathsf{\ while\ points\ of\ }\,\ns'\mathsf{\ satisfy\ }(\alpha).
\end{equation}
The first claim is obvious from (\ref{sng}), since $\,\sg\hh'$ 
consists of singular points of $\,Z$, cf. \cite[Remark 6.2(a)]{derdzinski}. 
As for the second one, its failure would---again by (\ref{sng})---amount 
to ($\beta$) for some points $\,y\in\ns'\nnh$, arbitrarily close to $\,x$. 
Combined with Theorem~\ref{zrset}, this would imply that 
$\,\tyn'\nh=\mathcal{H}_y^{\phantom i}\nh\cap\mathcal{H}_y^\perp\nnh$. 
Since $\,\tyn'\nh=\mathrm{Ker}\hskip1.7pt\nabla\nh v_y^{\phantom i}$ by 
(\ref{tyn}), both inclusions
\[
\mathcal{H}_y^{\phantom i}\nh\cap\mathcal{H}_y^\perp\subset\,
\mathcal{H}_y^{\phantom i}=\mathrm{Ker}\hskip1.7pt\nabla\nh v_y^{\phantom i}
\cap\hs\mathrm{Ker}\hskip2.2ptd\phi_y^{\phantom i}\hskip8pt\mathsf{and}
\hskip8pt
\mathcal{H}_y^{\phantom i}=\mathrm{Ker}\hskip1.7pt\nabla\nh v_y^{\phantom i}
\cap\hs\mathrm{Ker}\hskip2.2ptd\phi_y^{\phantom i}\hs
\subset\,\mathrm{Ker}\hskip1.7pt\nabla\nh v_y^{\phantom i}
\]
would be equalities. As $\,\phi(y)=0\,$ by (\ref{fez}), the second 
inclusion-turned-equality would give 
$\,\nabla\nh\phi_y^{\phantom i}\hs\in\hs\nabla\nh v_y^{\phantom i}(\tym)$, 
and so Theorem~\ref{esszr} would yield case 
($\alpha$) for $\,y\,$ rather than ($\beta$). The ensuing contradiction 
proves the second part of (\ref{alg}).

\ 

Let $\,\varPi_\alpha^{\phantom i}$ (or $\,\varPi_\beta^{\phantom i}$, or 
$\,\varPi_\gamma^{\phantom i}$) denote the subset of a given component 
$\,\ns\,$ of $\,Z\,$ formed by all points $\,x\in\ns\hs$ with 
($\alpha$) (or ($\beta$) or, respectively, ($\gamma$)). According to 
\cite[Remark 17.1]{derdzinski}, $\,\varPi_\alpha^{\phantom i}$ and 
$\,\varPi_\beta^{\phantom i}$ are relatively open in $\,\ns\hn$. So is, 
consequently, the set $\,\ns'\nh\cup\sg\hh'\nh
=\exp\hskip-2pt_x^{\phantom i}[\hs C\cap H\cap U\hs]\,$ appearing in 
(\ref{fez}) (by Theorem~\ref{zrset}), as well as the union 
$\,\varPi_\alpha^{\phantom i}\nh\cup\varPi_\gamma^{\phantom i}$ (in view of 
(\ref{alg})). Thus, due to connectedness of $\,\ns\hn$, either 
$\,\ns=\varPi_\beta^{\phantom i}$ (in which case $\,\ns\,$ is essential, and 
we denote it by the symbol $\,\es$), or 
$\,\ns=\varPi_\alpha^{\phantom i}\nh\cup\varPi_\gamma^{\phantom i}$ is a 
nonessential component (and we let $\,\sg\,$ stand for the set of its 
essential points, so that $\,\sg=\varPi_\gamma^{\phantom i}$). In other words, 
since $\,\varPi_\alpha^{\phantom i},\hs\varPi_\beta^{\phantom i}$ and 
$\,\varPi_\gamma^{\phantom i}$ are pairwise disjoint, we have
\begin{equation}\label{ese}
(*)\hskip9pt\es=\varPi_\beta^{\phantom i}\hskip8pt\mathsf{in\ case\ (a),}
\hskip20pt(**)\hskip9pt\sg
=\varPi_\gamma^{\phantom i}\hskip5pt\mathsf{and}\hskip5pt
\ns\smallsetminus\sg=\varPi_\alpha^{\phantom i}\hskip8pt\mathsf{in\ case\ (b).}
\end{equation}
Assertions (i) -- (iii), in both cases (a) and (b), are now immediate (with 
one exception): for any $\,x\in\es$, Theorem~\ref{zrset} and (\ref{alg}) imply 
that $\,\sg\cap U'\nh=\sg\hh'\nnh$, with $\,\sg\hh'$ as above and 
sufficiently small $\,\,U,U'\nnh$. The exception is the possibility, still to 
be excluded, that, in case (b), $\,\sg\,$ might have connected components of 
different dimensions.

\ 

Assuming now case (b), we obtain (iv) as an obvious consequence of (\ref{sng}) 
and (\ref{ese}\hs-$**$). Now (iv) and Remark~\ref{tlumb} yield the first part 
of (v), while (vi) and (vii) follow from (\ref{tyn}), (\ref{alg}) and 
(\ref{ese}\hs-$**$). 

\

To prove the remainder of (v), first note that the claim about the sign 
pattern is true locally: a local con\-for\-mal change of the metric allows us 
to treat $\,v\,$ as a Kil\-ling field and use the final clause of 
Theorem~\ref{zrklf}, which implies that the the tangent spaces of 
$\,\ns\smallsetminus\sg\,$ are invariant under parallel transports along 
$\,\ns\smallsetminus\sg$. The corresponding global claim could fail only if 
some connected component of $\,\sg\,$ would locally 
disconnect $\,\ns\nnh$, leading to different sign patterns on the resulting 
new components. This, however, cannot happen since, for any $\,x\in\sg$, any 
$\,\ve\in(0,\infty)$, and any null geodesic 
$\,(-\hh\ve,\ve)\ni t\mapsto y(t)\,$ with $\,y(0)=x\,$ which lies in 
$\,\ns\smallsetminus\sg\,$ except at $\,t=0$, the family of tangent spaces 
$\,T\hskip-2.3pt_{y(t)}^{\phantom i}(\ns\smallsetminus\sg)$, for $\,t\ne0$, is 
parallel along the geodesic. Namely, whenever 
$\,t\mapsto w(t)\in T\hskip-2.3pt_{y(t)}\hskip-.9ptM\,$ is a parallel vector 
field and $\,g(\dot y,w)=0$, relations (\ref{nyw}) form a system of 
first-or\-der linear homogeneous ordinary differential equations with the 
unknowns $\,\nabla\nh_{\!w}^{\phantom i}\hn v\,$ and $\,g(w,\nabla\nh\phi)$. 
Therefore, if we choose a parallel field $\,w\,$ satisfying at some fixed 
$\,t\ne0\,$ the condition 
$\,w(t)\in T\hskip-2.3pt_{y(t)}(\ns\smallsetminus\sg)\,$ (so that, by (vi) and 
(\ref{fez}) -- (\ref{ese}), $\,\nabla\nh_{\!w}^{\phantom i}\hn v\,$ and 
$\,g(w,\nabla\nh\phi)\,$ both vanish at $\,t$), uniqueness of solutions gives 
$\,\nabla\nh_{\!w}^{\phantom i}\hn v=0\,$ and $\,g(w,\nabla\nh\phi)=0\,$ at 
every  $\,t$. Now, by (vi), 
$\,w(t)\in T\hskip-2.3pt_{y(t)}(\ns\smallsetminus\sg)\,$ for all $\,t\ne0$, as 
required. 

\ 

By (vii) the limit as $\,t\to0\,$ of the above parallel family 
$\,t\mapsto T\hskip-2.3pt_{y(t)}^{\phantom i}(\ns\smallsetminus\sg)\,$ is 
$\,u^\perp\hskip-1.2pt\cap H$, where $\,u=\dot y(0)\in\txm$, so that 
$\,u\in(C\cap H)\smallsetminus H^\perp\nnh$. Clearly, 
$\,\dim\hs(\ns\smallsetminus\sg)=\dim\hs(u^\perp\hskip-1.2pt\cap H)$. Letting 
$\,\sg\,$ temporarily stand for the connected component of $\,\sg\,$ which 
contains $\,x$, we have, by (ii), $\,\dim\sg=\dim\hh(H\cap H^\perp)$. Note 
that $\,H\cap H^\perp$ is the null\-space of the symmetric bilinear form 
$\,\langle\,,\rangle\,$ in $\,H\,$ obtained by restricting the metric $\,g$, 
and $\,r\,$ is the rank of the restriction of $\,\langle\,,\rangle\,$ to 
$\,u^\perp\hskip-1.2pt\cap H$. Since 
$\,u\in(C\cap H)\smallsetminus H^\perp\nnh$, the sign pattern of the latter 
restriction arises from that of $\,\langle\,,\rangle\,$ in $\,H\,$ by 
replacing a plus-minus pair with a zero. Consequently, $\,\langle\,,\rangle\,$ 
has the rank $\,r+2$, and $\,\dim\hs(\ns\smallsetminus\sg)\hs-\hs\dim\hs\sg
=\dim\hs(u^\perp\hskip-1.2pt\cap H)-\dim\hh(H\cap H^\perp)
=(\dim H-1)-[\hs\dim H-(r+2)]=r+1$. Now the dimension formula in (v) follows. 
This in turn shows that all connected components of $\,\sg\,$ have the same 
dimension, completing the proof.
}
\end{proof}

\section{Induced structures on $\,\sg\,$ and 
$\,\ns\smallsetminus\sg$}\label{is}
Again, let $\,\es\,$ now be either an essential component, or the set of 
essential points in a nonessential component $\,\ns\hs$ of the zero set 
$\,Z\,$ of a con\-for\-mal vector field $\,v\,$ on a 
pseu\-\hbox{do\hs-}\hskip0ptRiem\-ann\-i\-an manifold $\,(M,g)\,$ of dimension 
$\,n\ge3$.

\ 

Both $\,\es\,$ and $\,\ns\smallsetminus\sg\,$ carry geometric structures 
naturally induced by the underlying con\-for\-mal structure of $\,(M,g)$. 
Fixing our metric $\,g\,$ within the con\-for\-mal structure allows us in turn 
to represent the induced structures by more concrete geometric objects, as 
explained below.

\ 

First, according to Theorem~\ref{essen}(v), $\,g\,$ (or, the con\-for\-mal 
structure), restricted to $\,\ns\smallsetminus\sg$, is a symmetric 
$\,2$-tensor field having the same sign pattern at all points (or, 
respectively, a class of such tensor fields, arising from one another via 
multiplications by functions without zeros). We refer to it as the {\em 
pos\-si\-bly-de\-gen\-er\-ate metric\/} (or, {\em 
pos\-si\-bly-de\-gen\-er\-ate con\-for\-mal structure\/}) of 
$\,\ns\smallsetminus\sg$. If $\,\sg\subset\ns\,$ is nonempty, the 
met\-ric/struc\-ture must actually be degenerate due to the inequality in 
Theorem~\ref{essen}(v), which gives $\,r<\dim\hs(\ns\smallsetminus\sg)$, 
and also shows that this is the zero met\-ric/struc\-ture (with $\,r=0$) only 
in the case where $\,\dim\hs\sg\hs=\dim\hs(\ns\smallsetminus\sg)\hs-\hs1$.

\ 

A further natural structure on $\,\ns\smallsetminus\sg\,$ is the {\em 
null\-space distribution\/} $\,\mathcal{P}\,$ of the restriction of the metric 
$\,g\,$ (or con\-for\-mal structure) to $\,\ns\smallsetminus\sg$. The 
dimension of $\,\mathcal{P}\,$ is positive when $\,\sg\subset\ns\,$ is 
nonempty: as we just saw, the restricted metric is degenerate. Instead of the 
argument in the last paragraph, we can also derive its degeneracy from the 
fact that, by the Gauss lemma, short null geodesic segments emanating from 
$\,\sg\,$ into $\,\ns\smallsetminus\sg\,$ are all tangent to $\,\mathcal{P}$. 
(The Gauss lemma and its standard proof in the Riemannian case 
\cite[Lemma 10.5]{milnor} remain valid for indefinite metrics.)

\ 

From now on $\,\sg\,$ is assumed nonempty. In view of Theorem~\ref{essen}(i), 
$\,g\,$ gives rise to an obvious tor\-sion-free connection $\,\mathrm{D}\,$ on 
$\,\es$, while the con\-for\-mal structure of $\,g\,$ induces on $\,\es\,$ a 
natural {\em projective structure}, that is, a class of tor\-sion-free 
connections having the same family of non\-pa\-ram\-e\-triz\-ed geodesics. 
(See the third paragraph of Section~\ref{pr}.)

\ 

In addition, $\,g\,$ naturally leads to a $\,1$-form $\,\xi\,$ on $\,\es$. 
(Using the con\-for\-mal structure instead of $\,g$, we obtain a $\,1$-form 
$\,\xi\,$ defined only up to multiplications by functions without zeros.) 
To describe $\,\xi$, we consider two cases, noting that 
$\,\phi=(2/n)\hskip2pt\mathrm{div}\hskip2ptv\,$ is constant on $\,\es\,$ 
and, in fact, on every component of $\,Z$, cf.\ the final clause of 
Theorem~\ref{zrset}. Specifically, if $\,\phi=0\,$ on $\,\es$, then 
$\,\es\ni x\,\mapsto\,\mathcal{H}_x^{\phantom i}
=\hs\mathrm{Ker}\hskip1.7pt\nabla\nh v_x^{\phantom i}
\cap\hs\mathrm{Ker}\hskip2.2ptd\phi_x^{\phantom i}$ is, in both cases, 
a parallel subbundle of $\,T\hskip-3.2pt_\sg^{\phantom i}M$ contained in 
$\,\mathrm{Ker}\hskip1.7pt\nabla\nh v\,$ as a codimension-one subbundle 
\cite[Lemma 13.1(b),\hs(d)]{derdzinski}, and 
we set $\,\xi=g(u,\,\cdot\,)$, on $\,\es$, for any section $\,u\,$ of 
$\,\mathrm{Ker}\hskip1.7pt\nabla\nh v$ over $\,\es\,$ with 
$\,g(u,\nabla\nh\phi)=1$. If $\,\phi\ne0\,$ on $\,\es$, we declare that 
$\,\xi=0$.
\begin{proposition}\label{pties}Under the assumptions made at the beginning of 
this section, for\/ $\,\mathcal{P},\hs\mathrm{D}\,$ and\/ $\,\xi\,$ defined 
above,
\begin{enumerate}
  \def\theenumi{{\rm\roman{enumi}}}
\item[{\rm(i)}] $\mathcal{P}\,$ is integrable and its leaves are  null totally 
geodesic sub\-man\-i\-folds of\/ $\,(M,g)$,
\item[{\rm(ii)}] if\/ $\,\varGamma\subset\sg\,$ is a geodesic segment and\/ 
$\,T\hskip-2.3pt_x^{\phantom i}\varGamma
\subset\mathrm{Ker}\hskip2.7pt\xi_{\hs x}^{\phantom i}$ 
for some\/ $\,x\in\varGamma\nnh$, then\/ 
$\,T\hskip-2.3pt_x^{\phantom i}\varGamma
\subset\mathrm{Ker}\hskip2.7pt\xi_{\hs x}^{\phantom i}$ for every\/ 
$\,x\in\varGamma\nnh$,
\item[{\rm(iii)}] in the open subset\/ $\,\es\hh'\nh\subset\es\,$ on which\/ 
$\,\xi\ne0$,
\begin{equation}\label{sym}
\mathrm{sym}\hskip1.7pt\mathrm{D}\hs\xi\,\,=\,\,\mu\odot\hs\xi\hs,\hskip10pt
\text{that\ is,}\hskip8pt\xi_{\hs j,\hs k}^{\phantom i}
+\xi_{\hs k,\hskip1.3ptj}^{\phantom i}
=\mu_j^{\phantom i}\xi_{\hs k}^{\phantom i}
+\mu_k^{\phantom i}\xi_{\hs j}^{\phantom i}\hskip10pt\text{for\ 
some\ }\,1\hyp\text{form\ }\,\,\mu\,\hs\text{\ on\ }\,\hs\sg\hh',
\end{equation}
\item[{\rm(iv)}] $\xi\,$ has the following unique continuation 
property\/{\sf:} if\/ $\,\xi=0\,$ at all points of some co\-di\-men\-sion-one 
connected submanifold\/ $\,\Delta\,$ of\/ $\,\es$, then\/ $\,\xi=0\,$ 
everywhere in the connected component of\/ $\,\es\,$ containing\/ 
$\,\Delta$.
\end{enumerate}
\end{proposition}
\begin{proof}
{\sf 
For any sections $\,w,w'$ of $\,\mathcal{P}\,$ and any curve 
$\,t\mapsto y(t)\,$ in the totally umbilical sub\-man\-i\-fold 
$\,K=\ns\smallsetminus\sg$, (\ref{bxw}) gives 
$\,\pi\nabla_{\!\dot x}^{\phantom i}w=0$, that is, 
$\,\nabla_{\!\dot x}^{\phantom i}w\,$ is tangent to $\,K$. Hence so is 
$\,\nabla\nh_{\!w}^{\phantom i}w'$ and, for any vector field $\,u\,$ 
tangent to $\,K\,$ we have 
$\,g(\nabla\nh_{\!w}^{\phantom i}w',u)
=-\hh g(\nabla\nh_{\!w}^{\phantom i}u,w')=0$, as one sees applying 
(\ref{bxw}), this time, to $\,w'$ instead of $\,w\,$ and an integral 
curve $\,t\mapsto y(t)\,$ of $\,w$. Thus, $\,\nabla\nh_{\!w}^{\phantom i}w'$ 
is a section of $\,\mathcal{P}\nnh$, and (i) follows.

\ 

To prove (ii) -- (iv), we may assume that $\,\phi=0\,$ on $\,\es$. For 
$\,\varGamma\,$ and $\,x\,$ as in (ii), let 
$\,t\mapsto w(t)\in T\hskip-2.3pt_{y(t)}\hskip-.9ptM\,$ be a parallel vector 
field along a geodesic parametrization of $\,t\mapsto y(t)\,$ of 
$\,\varGamma\,$ such that $\,y(0)=x\,$ and 
$\,\nabla_{\!w(0)}^{\phantom i}v=\dot y(0)$. 
(That $\,\dot y(0)\in\nabla\nh v_x^{\phantom i}(\txm)\,$ is clear: as 
$\,\dot y(0)\,$ lies in $\,T\hskip-2.3pt_x^{\phantom i}\varGamma
\subset\mathrm{Ker}\hskip2.7pt\xi_{\hs x}^{\phantom i}
\subset T\hskip-2.3pt_x\es$, it is orthogonal not just to 
$\,\mathcal{H}_x^{\phantom i}
=\mathrm{Ker}\hskip1.7pt\nabla\nh v_x^{\phantom i}
\cap\hs\mathrm{Ker}\hskip2.2ptd\phi_x^{\phantom i}$, cf.\ 
Theorem~\ref{essen}(ii), but also to the whole space 
$\,\mathrm{Ker}\hskip1.7pt\nabla\nh v_x^{\phantom i}$, while 
$\,\nabla\nh v_x^{\phantom i}(\txm)
=[\hs\mathrm{Ker}\hskip1.7pt\nabla\nh v_x^{\phantom i}]^\perp$ as 
$\,\nabla\nh v_x^{\phantom i}$ is skew-adjoint by (\ref{tnv}) with 
$\,\phi(x)=0$.) Choosing a function $\,t\mapsto\kappa(t)\,$ with 
$\,2\hh\dot\kappa=g(w,\nabla\nh\phi)\,$ and $\,\kappa(0)=1$, then integrating 
(\ref{nyw}\hs-a), we get 
$\,\nabla_{\!w(t)}^{\phantom i}v=\kappa(t)\dot y(t)$, so that 
$\,\dot y(t)\in\nabla\nh v_{y(t)}^{\phantom i}(\txm)
=[\hs\mathrm{Ker}\hskip1.7pt\nabla\nh v_{y(t)}^{\phantom i}]^\perp$ for 
all $\,t\,$ near $\,0$, which yields (ii).

\ 

Assertion (iii) is in turn a consequence of (ii): if $\,x\in\es\,$ and 
$\,w\in\mathrm{Ker}\hskip2.7pt\xi_{\hs x}^{\phantom i}$, setting 
$\,y(t)=\exp\hskip-2pt_x^{\phantom i}\hs tw\,$ and differentiating the 
resulting equality $\,\xi(\dot y)=0$, we obtain 
$\,[\nabla_{\!\dot y}^{\phantom i}\xi\hh](\dot y)=0$, so that 
$\,[\nabla_{\!t y}^{\phantom i}\xi\hh](\dot y)=0$, so that 
$\,[\nabla\nh_{\!w}^{\phantom i}\xi\hh](w)=0$. Thus, 
$\,\mathrm{sym}\hskip1.7pt\mathrm{D}\hs\xi_x^{\phantom i}$ treated as a 
polynomial function on $\,T\hskip-2.3pt_x^{\phantom i}\es\,$ vanishes on 
the zero set of the linear function $\,\xi_x^{\phantom i}$, which is 
well-known to imply divisibility of the former by the latter, 
cf.\ \cite[Lemma 17.1(i)]{derdzinski-maschler}, thus proving (iii).

\ 

Finally, (iv) follows from (ii) since under 
the hypothesis of (iv), $\,\xi\,$ must vanish on an open set containing 
$\,\Delta$, namely, the set of points at which an open set of tangent 
directions is realized by geodesics intersecting $\,\Delta$.
}
\end{proof}
Note that condition (\ref{sym}) involves $\,\mathrm{D}\,$ only through its 
underlying projective structure, and remain valid after $\,\xi\,$ has been 
multiplied by a function without zeros.

\section{One\hh-\nnh jets of $\,v\,$ along components of $\,Z$}\label{oj}
As before, $\,Z\,$ stands for the zero set of a con\-for\-mal vector field 
$\,v\,$ on a pseu\-\hbox{do\hs-}\hskip0ptRiem\-ann\-i\-an manifold $\,(M,g)\,$ 
of dimension $\,n\ge3$. For $\,x\in Z$, the en\-do\-mor\-phism 
$\,\nabla\nh v_x^{\phantom i}$ of $\,\txm$, obviously independent of the 
choice of the connection $\,\nabla\nnh$, is also known as the {\em linear 
part\/} (or {\em Jacobian}, or {\em derivative}, or {\em differential\/}) of 
$\,v\,$ at the zero $\,x$. At the same time, $\,\nabla\nh v_x^{\phantom i}$ is 
the infinitesimal generator of the local flow of $\,v\,$ acting in $\,\txm$.

\ 

Since $\,v_x^{\phantom i}=0$, we may also identify 
$\,\nabla\nh v_x^{\phantom i}$ with the $\,1$-jet of $\,v\,$ at $\,x$.

\ 

Given $\,x,y\in Z$, we say that the $\,1$-jets of $\,v\,$ at $\,x\,$ and 
$\,y\,$ are {\em con\-for\-mal\-ly equivalent\/} if, for some 
vertical\hh-arrow con\-for\-mal isomorphism $\,\txm\to\tym$, the following 
diagram commutes:
\begin{equation}\label{dia}
\begin{matrix}
\txm&\mathop{-\!\!\!-\!\!\!-\!\!\!-\!\!\!-\!\!\!-\!\!\!\rightarrow}
\limits^{\nabla\nh v_x^{\phantom i}}&\txm\phantom{_|}\\
{}\vbox{\hbox{\Bigg\downarrow}\vskip-12pt}&
&\vbox{\hbox{\Bigg\downarrow}\vskip-12pt}\\
\tym&\mathop{-\!\!\!-\!\!\!-\!\!\!-\!\!\!-\!\!\!-\!\!\!\rightarrow}
\limits^{\nabla\nh v_y^{\phantom i}}&\tym
\end{matrix}
\end{equation}
By {\em con\-for\-mal isomorphisms\/} we mean here nonzero scalar multiples of 
linear isometries.
\begin{proposition}\label{onjts}Under the assumptions of 
Theorem\/~{\rm\ref{essen}}, with\/ $\,\xi\,$ defined in 
Section\/~{\rm\ref{is}},
\begin{enumerate}
  \def\theenumi{{\rm\roman{enumi}}}
\item[{\rm(i)}] in case\/ {\rm(b)} of Theorem\/~{\rm\ref{essen}}, for any 
connected component\/ $\,\ns'$ of\/ $\,\ns\smallsetminus\sg$, the\/ $\,1$-jets 
of\/ $\,v\,$ at all points of\/ $\,\ns'$ are con\-for\-mal\-ly equivalent to 
one another, but not to the\/ $\,1$-jet of\/ $\,v\,$ at any point of\/ $\,\sg$,
\item[{\rm(ii)}] in both cases\/ {\rm(a)} -- {\rm(b)} of 
Theorem\/~{\rm\ref{essen}}, if\/ $\,\xi\,$ is not identically zero on a 
connected component\/ $\,\es\hh'$ of\/ $\,\es$, then the\/ $\,1$-jets of\/ 
$\,v\,$ at any two points of\/ $\,\es\hh'$ are con\-for\-mal\-ly equivalent.
\end{enumerate}
\end{proposition}
\begin{proof}
{\sf 
Since some local con\-for\-mal change of the metric near any $\,y\in\ns'$ 
turns $\,v\,$ into a Kil\-ling field, (i) follows: $\,\nabla\nh v\,$ then 
becomes parallel along a neighborhood of $\,y\,$ in $\,\ns'\nnh$. The claim 
about $\,\sg\,$ is obvious from Theorem~\ref{essen}(vi).

\ 

In (ii), the definition of $\,\xi\,$ (see Section~\ref{is}) implies that 
$\,\phi=0\,$ on $\,\es\hh'\nnh$. Thus, by Theorem~\ref{esszr}, 
$\,\nabla\nh\phi_x^{\phantom i}\hs\notin\hs\nabla\nh v_x^{\phantom i}(\txm)\,$ 
at every point $\,x\in\es\hh'\nnh$, and so 
$\,\mathcal{H}_x^\perp\nh=\nabla\nh v_x^{\phantom i}(\txm)
\oplus\bbR\nabla\nh\phi_x^{\phantom i}$. (Note that 
$\,\nabla\nh v_x^{\phantom i}(\txm)
=[\hs\mathrm{Ker}\hskip1.7pt\nabla\nh v_x^{\phantom i}]^\perp$ by 
(\ref{tnv}) with $\,\phi(x)=0$, while (\ref{hxe}) gives 
$\,\nabla\nh\phi_x^{\phantom i}\in\mathcal{H}_x^\perp\nnh$.) Hence, in view of 
Theorem~\ref{essen}(ii), 
$\,T\hskip-2.3pt_x^{\phantom i}\es\hh'\nh
\subset\nabla\nh v_x^{\phantom i}(\txm)
\oplus\bbR\nabla\nh\phi_x^{\phantom i}$, while the vectors tangent to 
$\,\es\hh'$ at $\,x\,$ and lying in the summand 
$\,\nabla\nh v_x^{\phantom i}(\txm)\,$ form precisely the sub\-space 
$\,\mathrm{Ker}\hskip2.7pt\xi_{\hs x}^{\phantom i}$, which has 
co\-di\-men\-sion one in $\,T\hskip-2.3pt_x^{\phantom i}\es\hh'$ for all 
points $\,x\,$ of a dense open subset of $\,\es\hh'$ 
(Proposition~\ref{pties}(iv)). We will now show that the con\-for\-mal 
equivalence type of the $\,1$-jets of $\,v\,$ is constant along any 
geodesic segment $\,\varGamma\,$ in $\,\es\hh'$ with a parametrization 
$\,t\mapsto y(t)$ satisfying the condition 
$\,\dot y(t)\notin\hs\nabla\nh v_x^{\phantom i}(\txm)\,$ at each 
$\,x=y(t)$. (As any two points of $\,\es\hh'$ can be joined by piecewise 
smooth curves made up from such geodesic segments, due to the denseness and 
openness property just mentioned, (ii) will then clearly follow.)

\ 

Specifically, our assumption about $\,\dot y(t)\,$ yields 
$\,\nabla\nh\phi=\rho\hs\dot y+\nabla\nh_{\!w}^{\phantom i}\hn v\,$ for some 
function $\,t\mapsto\rho(t)\,$ and a vector field 
$\,t\mapsto w(t)\in T\hskip-2.3pt_{y(t)}\hskip-.9ptM\,$ along the geodesic; 
since $\,\mathrm{rank}\hskip1.7pt\nabla\nh v\,$ is constant on $\,\es\hh'$ 
by \cite[Lemma 13.1(d)]{derdzinski}, $\,w\,$ may be chosen differentiable.) 
From (\ref{nwn}) we now obtain $\,2\hh\nabla_{\!\dot y}^{\phantom i}\nabla\nh v
=g(\nabla\nh\phi,\,\cdot\,)\nnh\otimes\nnh\dot y
-g(\dot y,\,\cdot\,)\nnh\otimes\nnh\nabla\nh\phi$, and one easily verifies 
that $\,\nabla\nh v\,$ is $\,\mathrm{D}$-parallel for the new metric 
connection $\,\mathrm{D}\,$ in $\,T\hskip-2.4pt_\varGamma^{\phantom i}M\,$  
given by $\,2\hs\mathrm{D}\nh_{\dot y}=2\hh\nabla_{\!\dot y}
+g(w,\,\cdot\,)\nnh\otimes\nnh\dot y-g(\dot y,\,\cdot\,)\nnh\otimes\nnh w$.
}
\end{proof}

\section{The associated quintuples}\label{aq}
The symbol $\,[\ip]\,$ stands for the {\em homothety class\/} of a 
pseu\-\hbox{do\hs-}\hskip0ptEuclid\-e\-an inner product $\,\ip\,$ on a 
fi\-\hbox{nite\hh-}\hskip0ptdi\-men\-sion\-al 
vector space $\,\pc\nnh$, that is, the set of all nonzero scalar multiples of 
$\,\ip$. The underlying con\-for\-mal structure $\,[\hs g\hh]\,$ of 
a pseu\-\hbox{do\hs-}\hskip0ptRiem\-ann\-i\-an manifold $\,(M,g)\,$ may thus 
be identified with the assignment $\,M\ni x\mapsto[\hs g_x]$. Let us consider 
quintuples
\begin{equation}\label{qui}
(\pc\nnh,\,[\ip],\,B,\,\lambda,\,\dl)
\end{equation}
formed by a 
pseu\-\hbox{do\hs-}\hskip0ptEuclid\-e\-an vector space $\,\pc\nnh$, the 
homothety class of its inner product $\,\ip$, a skew-ad\-joint 
en\-do\-mor\-phism $\,B\in\mathfrak{so}\hh(\pc)$, a real number $\,\lambda$, 
and a linear functional 
$\,\dl\in[\hh\mathrm{Ker}\hskip1.7pt(B+\lambda)]^*$ on the 
sub\-space $\,\mathrm{Ker}\hskip1.7pt(B+\lambda)\,$ of $\,\pc$ (which, if 
nontrivial, is the eigen\-space of $\,B\,$ for the eigen\-value $\,-\lambda$).

\ 

We call (\ref{qui}) {\em algebraically equivalent\/} to another such 
quintuple $\,(\pc\hh'\nnh,\,[\ip\hh'\hh],\,B\hh'\nnh,\,\lambda'\nnh,
\,\dl\hh'\hh)\,$ if $\,\lambda'\nh=\lambda\,$ and some linear isomorphism 
$\,\pc\to\pc\hh'$ sends $\,[\ip],\,B,\dl\,$ to 
$\,[\ip\hh'\hh],\,B\hh'$ and $\,\dl\hh'\nnh$.

\ 

Examples of quintuples (\ref{qui}) arise as follows. Given a con\-for\-mal 
vector field $\,v\,$ on a pseu\-\hbox{do\hs-}\hskip0ptRiem\-ann\-i\-an 
manifold $\,(M,g)$ and a point $\,x\in M\,$ at which $\,v_x\nh=0$, we define 
the {\em quintuple associated with\/} $\,v\hs$ {\em and\/} $\,x\,$ to be 
$\,(\pc\nnh,\,[\ip],\,B,\,\lambda,\,\dl)
=(\txm\nh,\,[\hs g_x],\,A_x,\,\phi(x),\dl)$, where $\,A\,$ and $\,\phi\,$ are 
determined by $\,v\,$ as in (\ref{tnv}), and $\,\dl\,$ is the restriction of 
$\,d\phi_x$ to the subspace 
$\,\mathrm{Ker}\hskip1.7pt(B+\lambda)
=\mathrm{Ker}\hskip1.7pt\nabla\nh v_x$. In other words, $\,B\,$ equals twice 
the skew-ad\-joint part of $\,\nabla\nh v_x:\txm\to\txm\,$ (the value at 
$\,x\,$ of the morphism $\,\nabla\nh v:\tm\to\tm$), and $\,\lambda\,$ is 
$\,2/n\,$ times $\,\mathrm{tr}\hskip1.7pt\nabla\nh v_x$, where $\,n=\dim M\nh$.

\ 

The associated quintuple $\,(\pc\nnh,\,[\ip],\,B,\,\lambda,\,\dl)\,$ depends, 
besides $\,v\,$ and $\,x$, only on the underlying con\-for\-mal structure 
$\,[\hs g\hh]$, rather than the metric $\,g$. This is obvious for 
$\,\pc,[\ip],B\,$ and 
$\,\lambda=(2/n)\hskip2pt\mathrm{tr}\hskip1.7pt\nabla\nh v_x$, cf.\ the 
beginning of Section~\ref{oj}. Similarly, as 
$\,\mathrm{Ker}\hskip1.7pt(B+\lambda)
=\mathrm{Ker}\hskip1.7pt\nabla\nh v_x$, the last line in Remark~\ref{phcnf} 
yields the claim about $\,\dl$.

\section{Con\-for\-mal equivalence of two\hs-\nh jets}\label{ce}
Let $\,v\,$ and $\,w\,$ be con\-for\-mal vector fields on 
pseu\-\hbox{do\hs-}\hskip0ptRiem\-ann\-i\-an manifolds $\,(M,g)\,$ and, 
respectively, $\,(N\nh,h)$, such that $\,v\,$ vanishes at a point 
$\,x\in M\nh$, and $\,w\,$ at $\,y\in N\nh$. We say that the $\,2$-jet 
of $\,v\,$ at $\,x\,$ is {\em con\-for\-mal\-ly equivalent\/} to the $\,2$-jet 
of $\,w\,$ at $\,y\,$ if some dif\-feo\-mor\-phism $\,F\,$ between a 
neighborhood $\,\,U\,$ of $\,x\,$ in $\,M\,$ and one of $\,y\,$ in 
$\,N\nh$, with $\,F(x)=y$, sends the former $\,2$-jet to the latter, 
while, at the same time, for some function $\,\tau:U\to\bbR$, the metrics 
$\,F\hh^*\nnh h\,$ and $\,e^\tau\hskip-2ptg$ have the same $\,1$-jet at $\,x$.

\ 

As $\,v\,$ and $\,w\,$ vanish at $\,x\,$ and $\,y$, the above condition 
on $\,F\,$ involves $\,F\,$ only through its $\,2$-jet at $\,x$.
\begin{lemma}\label{equiv}For\/ $\,M,g,v,x\,$ and\/ $\,N\nh,h,w,y\,$ as in the 
last two paragraphs, the\/ $\,2$-jets of\/ $\,v\,$ at\/ $\,x\,$ and of\/ 
$\,w\,$ at\/ $\,y\,$ are con\-for\-mal\-ly equivalent if and only if the 
quintuple\/ 
$\,(\pc\nnh,\,[\ip],\,B,\,\lambda,\,\dl)\,$ associated with\/ 
$\,v\,$ and\/ $\,x\,$ is algebraically equivalent, in the sense of 
Section\/~{\rm\ref{aq}}, to the analogous quintuple\/ 
$\,(\pc\hh'\nnh,\,[\ip\hh'\hh],\,B\hh'\nnh,\,\lambda'\nnh,
\,\dl\hh'\hh)\,$ for\/ $\,w\,$ and\/ $\,y$.
\end{lemma}
\begin{proof}
{\sf 
The `only if' part of our claim is obvious from functoriality of 
the associated quintuple. To prove the `if' part, we fix local coordinates 
$\,x^{\hs j}$ for $\,M\,$ at $\,x\,$ and $\,y^{\hh a}$ for $\,N\hs$ at $\,y\,$ 
such that the corresponding Christoffel symbols of $\,g$, or $\,h\,$ vanish at 
$\,x$, or $\,y$. We also set 
$\,v_j^{\hs k}=\partial_j^{\phantom i}v^{\hs k}\nnh$, 
$\,F_{\hskip-2ptj}^{\hs a}=\partial_j^{\phantom i}F^{\hs a}\nnh$, 
$\,F_{\hskip-2ptjk}^{\hs a}
=\partial_j^{\phantom i}\partial_k^{\phantom i}F^{\hs a}\nnh$, 
$\,\tau_j^{\phantom i}=\partial_j^{\phantom i}\tau$, 
$\,w_c^{\hh a}=\partial_c^{\phantom i}w^{\hh a}\nnh$, where 
all the partial derivatives stand for their values at $\,x\,$ or $\,y\,$ (and 
those involving $\,F\,$ or $\,\tau\,$ are treated as unknowns).

It now suffices to show that, if 
$\,(\pc\nnh,\,[\ip],\,B,\,\lambda,\,\dl)\,$ and 
$\,(\pc\hh'\nnh,\,[\ip\hh'\hh],\,B\hh'\nnh,\,\lambda'\nnh,
\,\dl\hh'\hh)\,$ are equivalent, the system
\begin{equation}\label{sys}
\begin{array}{rl}
\mathrm{i)}&\hskip7pt
w_c^{\hh a}F_{\hskip-2ptj}^{\hs c}=F_{\hskip-1.4ptk}^{\hs a}v_j^{\hs k}\hh,\\
\mathrm{ii)}&\hskip7ptF_{\hskip-1.4ptl}^{\hs a}
\,\partial_j^{\phantom i}\partial_k^{\phantom i}v^{\hs l}\hs
+\,F_{\hskip-2ptjl}^{\hs a}v_k^{\hs l}
+\,F_{\hskip-1.4ptkl}^{\hs a}v_j^{\hs l}\hs
=\,F_{\hskip-2ptj}^{\hs b}F_{\hskip-1.4ptk}^{\hs c}
\,\partial_{\hs b}^{\phantom i}\partial_c^{\phantom i}w^{\hh a}\hs
+\,F_{\hskip-2ptjk}^{\hs c}w_c^{\hh a}\hh,\\
\mathrm{iii)}&\hskip7pt
h_{ac}^{\phantom.}F_{\hskip-2ptj}^{\hs a}F_{\hskip-1.4ptk}^{\hs c}
=e^\tau\hskip-2ptg_{jk}^{\phantom.}\hh,
\hskip16pt\mathrm{iv)}\hskip12pth_{ac}^{\phantom.}(F_{\hskip-2ptj}^{\hs a}
F_{\hskip-1.4ptkl}^{\hs c}
+F_{\hskip-2ptjl}^{\hs a}F_{\hskip-1.4ptk}^{\hs c})
=e^\tau\nnh\tau_l^{\phantom i}g_{jk}^{\phantom.}\hh,
\end{array}
\end{equation}
where the values of $\,g_{jk}^{\phantom.},
\hs\partial_j^{\phantom i}\partial_k^{\phantom i}v^{\hs l},
\hs h_{ac}^{\phantom.}$ and 
$\,\partial_{\hs b}^{\phantom i}\partial_c^{\phantom i}w^{\hh a}$ are taken at 
$\,x\,$ or $\,y$, has a solution consisting of a real number $\,\tau$ and some 
quantities $\,F_{\hskip-2ptj}^{\hs a}$, $\,F_{\hskip-2ptjk}^{\hs a}$, 
$\,\tau_j^{\phantom i}$ with 
$\,F_{\hskip-2ptjk}^{\hs a}=F_{\hskip-1.4ptkj}^{\hs a}$.

\ 

As the first part of such a solution we choose a matrix 
$\,F_{\hskip-2ptj}^{\hs a}$ which, when treated as a linear isomorphism 
$\,\txm\to\tyn$ (with the aid of our fixed coordinates $\,x^{\hs j}$ and 
$\,y^{\hh a}$), realizes the equivalence of the two quintuples. As 
$\,\lambda'\nh=\lambda\,$ and the isomorphism in question sends 
$\,([\ip],\,B,\,\dl)\,$ to 
$\,([\ip\hh'\hh],\,B\hh'\nnh,\,\dl\hh'\hh)$, we now clearly have 
(\ref{sys}\hs-iii) for some $\,\tau\in\bbR$, (\ref{sys}\hs-i), and there 
exists a $\,1$-form $\,\sigma\in\taxm\,$ with 
$\,\phi_{\nh,\hs j}^{\phantom i}
-\psi_{\nh,\hs a}^{\phantom i}F_{\hskip-2ptj}^{\hs a}
=2\hh v_j^{\hs k}\sigma_k^{\phantom i}$, where $\,\phi\,$ is determined by 
$\,v\,$ as in (\ref{tnv}), $\,\psi\,$ is its analog for $\,w$, and the 
components of $\,d\phi\,$ and $\,d\hh\psi\,$ are evaluated at $\,x\,$ or 
$\,y$. (That $\,\sigma\,$ exists is obvious since our isomorphism sends 
$\,\dl\,$ to $\,\dl\hh'\nnh$, and so $\,u^{\hs j}(\phi_{\nh,\hs j}^{\phantom i}
-\psi_{\nh,\hs a}^{\phantom i}F_{\hskip-2ptj}^{\hs a}\hh)=0\,$ whenever 
$\,u\in\txm\,$ and $\,u^{\hs j}v_j^{\hs k}=0$.) It follows that
\begin{equation}\label{gjk}
\mathrm{a)}\hskip5pt
g^{\hs jk}F_{\hskip-2ptj}^{\hs a}F_{\hskip-1.4ptk}^{\hs c}
=e^\tau\nnh h^{ac}\hh,\hskip11pt\mathrm{b)}\hskip5pt
v_k^{\hs j}\sigma^{\hs k}\nnh=\phi\hs\sigma^{\hs j}
-\hh g^{jl}v_l^{\hs k}\sigma_k^{\phantom i}\hh,\hskip6pt\mathrm{where}\hskip6pt
\sigma^{\hs j}\nh=g^{jk}\sigma_k^{\phantom i}\hh.
\end{equation}
In fact, (\ref{sys}\hs-iii) states that the matrix 
$\,e^{-\tau\nnh/2}F_{\hskip-2ptj}^{\hs a}\nh$, as a linear isomorphism 
$\,\txm\to\tyn\nh$, sends the metric $\,g_x$ to $\,h_y$, and so the 
reciprocal metrics of $\,g_x$ and $\,h_y$ correspond to each other under the 
dual isomorphism $\,\tayn\to\taxm$. This amounts to (\ref{gjk}\hs-a), while 
(\ref{gjk}\hs-b) is a trivial consequence of (\ref{tnv}).

\ 

Finally, let us set 
$\,F_{\hskip-2ptjk}^{\hs a}=F_{\hskip-1.2ptl}^{\hs a}\sigma^{\hs l}\nnh
g_{jk}^{\phantom.}\nh-\sigma_j^{\phantom i}F_{\hskip-1.4ptk}^{\hs a}
-\sigma_k^{\phantom i}F_{\hskip-2ptj}^{\hs a}$ and 
$\,\tau_j^{\phantom i}=-2\hs\sigma_j^{\phantom i}$. Then (\ref{sys}\hs-iii) 
implies (\ref{sys}\hs-iv). Next, at $\,x\,$ and $\,y$, our choice of coordinates 
and the coordinate form of (\ref{nwn}) yield 
$\,\partial_j^{\phantom i}\partial_k^{\phantom i}v^{\hs l}\nh
=\phi_{\nh,\hh k}^{\phantom i}\delta_{\nh j}^{\hs l}\nh
-\phi^{\hs,\hh l}g_{jk}^{\phantom.}\nh+\phi_{\nh,\hs j}^{\phantom i}\delta_k^{\hs l}$ and, 
analogously, $\,\partial_{\hs b}^{\phantom i}\partial_c^{\phantom i}w^{\hh a}
\nh=\psi_{\nh,\hh c}^{\phantom i}\delta_{\nh b}^{\hs a}\nh
-\psi^{\hs,\hh a}g_{bc}^{\phantom j}\nh
+\psi_{\nh,\hs b}^{\phantom i}\delta_c^{\hs a}$. 
Now (\ref{sys}\hs-iii) follows from (\ref{gjk}\hs-a), as one sees replacing 
$\,\phi_{\nh,\hs j}^{\phantom i}$ with 
$\,\psi_{\nh,\hs a}^{\phantom i}F_{\hskip-2ptj}^{\hs a}
+2\hh v_j^{\hs k}\sigma_k^{\phantom i}$ and noting that, by (\ref{tnv}), 
$\,v_j^{\hs l}F_{\hskip-2ptp}^{\hs a}\sigma^{\hs p}\nnh g_{kl}^{\phantom.}
+v_k^{\hs l}F_{\hskip-2ptp}^{\hs a}\sigma^{\hs p}\nnh g_{jl}^{\phantom.}
=\phi F_{\hskip-2ptp}^{\hs a}\sigma^{\hs p}\nnh g_{jk}^{\phantom.}$, while 
(\ref{sys}\hs-i) and (\ref{gjk}\hs-b) give 
$\,w_c^{\hh a}F_{\hskip-1.2ptl}^{\hs c}\sigma^{\hs l}\nh
=F_{\hskip-2ptp}^{\hs a}v_l^{\hh p}\sigma^{\hs l}\nh
=\phi F_{\hskip-2ptp}^{\hs a}\sigma^{\hs p}\nh
-g^{pq}F_{\hskip-2ptp}^{\hs a}v_q^{\hh l}\sigma_{\hn l}^{\phantom i}$.
}
\end{proof}

\section{Two\hh-\nnh jets of $\,v\,$ along components of $\,Z$}\label{tj}
Proposition~\ref{onjts} remains true if the word $\,1$-jet(s) is replaced 
everywhere with $\,2$-jet(s).

\ 

For both assertions (i) and (ii), this is a direct consequence of 
Lemma~\ref{equiv}. Specifically, in the case of (i), the invariant $\,\dl\,$ 
vanishes at every point of $\,\ns\smallsetminus\sg\,$ by Theorem~\ref{esszr}, 
while the remaining four objects in the associated quintuple (\ref{qui}) 
represent the $\,1$-jet of $\,v\,$ at the point in question. The 
al\-ge\-bra\-ic-e\-quiv\-a\-lence type of the quintuple is thus locally 
constant on $\,\ns\smallsetminus\sg\,$ as a consequence of 
Proposition~\ref{onjts}.

\ 

Similarly, for the new metric connection $\,\mathrm{D}\,$ in 
$\,T\hskip-2.4pt_\varGamma^{\phantom i}M\,$ used to prove 
Proposition~\ref{onjts}(ii), the second formula in (\ref{nwn}) shows that the 
restriction of $\,d\phi\,$ to 
$\,\mathrm{Ker}\hskip1.7pt\nabla\nh v\,$ is $\,\mathrm{D}$-parallel as well, 
as long as one chooses $\,w\,$ with $\,g(w,\nabla\nh\phi)=0$. Since 
$\,\nabla\nh\phi_x^{\phantom i}\hs\notin\hs\nabla\nh v_x^{\phantom i}(\txm)
=[\hs\mathrm{Ker}\hskip1.7pt\nabla\nh v_x^{\phantom i}]^\perp\nnh$, whenever 
$\,x\in\es\hh'\nnh$, such a choice is always possible.

\ 

The following example shows that the assumption about $\,\xi\,$ in 
Proposition~\ref{onjts}(ii) cannot in general be removed.

\ 

On a pseu\-\hbox{do\hs-}\hskip0ptEuclid\-e\-an space 
$\,(V\nnh,\langle\,,\rangle)\,$ of dimension $\,n\,$ we may define a 
con\-for\-mal vector field $\,v\,$ by
\begin{equation}\label{vxe}
v_x\,=\,\hs w\,+\,Bx\,+\,cx\,+\,2\langle u,x\rangle\hh x\,
-\,\langle x,x\rangle\hh u\hs,
\end{equation}
using any fixed vectors $\,w,u\in V\nnh$, any skew-ad\-joint 
en\-do\-mor\-phism $\,B$, and any scalar $\,c\in\bbR$. Let us now choose 
$\,n$ to be even, $\,\langle\,,\rangle\,$ to have the neutral signature, 
$\,B\,$ with two null $\,n$-dimensional eigenspaces for the eigenvalues 
$\,c,-\hh c$, and $\,u\,$ which does not lie in the $\,-\hh c\,$ eigenspace, 
along with $\,w=0$. Then 
$\,\dim\hs\mathrm{Ker}\hskip1.7pt\nabla\nh v_x^{\phantom i}$ is easily 
vefified to decrease when one replaces $\,x=0\,$ by any nearby vector $\,x\,$ 
orthogonal to $\,u\,$ and lying in the $\,-\hh c\,$ eigenspace of $\,B$.

\section*{Acknowledgements}
This paper is an expanded version of the author's talk in the first week of 
the Workshop on Cartan Connections, Geometry of Homogeneous Spaces, and 
Dynamics, which was organized by Andreas \v Cap, Charles Frances and Karin 
Melnick, and held at the Erwin Schr\"odinger International Institute for 
Mathematical Physics (ESI), Vienna, in July 2011. An essential part of the 
work on the paper was done during that week. The author thanks the University 
of Vienna, the ESI, and the organizers of the workshop for their support, 
hospitality, and the opportunity to participate in many stimulating 
discussions.

}

\end{document}